\documentclass[12pt]{amsart}

\usepackage{pinlabel}
\usepackage{url}
\usepackage{graphicx}
\usepackage{amssymb, amscd}
\usepackage{xy}
\usepackage{bm}
\xyoption{all}
\usepackage{xcolor}
\usepackage{scalerel}

\usepackage[OT2,OT1]{fontenc}
\newcommand\cyr{%
 \renewcommand\rmdefault{wncyr}%
 \renewcommand\sfdefault{wncyss}%
 \renewcommand\encodingdefault{OT2}%
 \normalfont
\selectfont}
 \DeclareTextFontCommand{\textcyr}{\cyr}

\DeclareMathOperator{\GL}{GL}

\newcommand{\N}{\mathbb N}
\newcommand{\Z}{\mathbb Z}

\newcommand{\R}{\mathbb R}

\newcommand{\CP}{\mathbb{CP}}

\newcommand{\Ghat}{\widehat G}
\newcommand{\rhat}{\widehat r}

\newcommand{\wSig}{\widetilde\Sigma}

\def\spinc{\ifmmode{\textrm{Spin}^c}\else{$\textrm{Spin}^c$}\fi}

\renewcommand{\phi}{\varphi}

\newtheorem{theorem}{Theorem}[section]
\newtheorem{thm}{Theorem}

\newtheorem{lemma}[theorem]{Lemma}

\newtheorem{proposition}[theorem]{Proposition}
\newtheorem{corollary}[theorem]{Corollary}
\theoremstyle{definition}

\newcommand{\veps}{\varepsilon}

\numberwithin{equation}{section}

\title{On the Thom conjecture in $\CP^3$}

\author[Daniel Ruberman]{Daniel Ruberman}
\address{Department of Mathematics, MS 050 \\ Brandeis University \\ Waltham, MA 02454}
\email{ruberman@brandeis.edu}

\author[Marko Slapar]{Marko Slapar}
\address{Faculty of Education\\University of Ljubljana\\Kardeljeva plo\v s\v cad 16\\1000
  Ljubljana, Slovenia\\ \newline
Faculty of Mathematics and Physics\\ University of Ljubljana \\ Jadranska 19\\1000 Ljubljana, Slovenia\\ and\newline
  Institute of Mathematics, Physics and Mechanics\\Jadranska
  19\\1000 Ljubljana, Slovenia}
\email{marko.slapar@pef.uni-lj.si}

\author[Sa\v{s}o Strle]{Sa\v{s}o Strle}
\address{Faculty of Mathematics and Physics, University of Ljubljana, Jadranska 19, 1000 Ljubljana, Slovenia, and\newline
  Institute of Mathematics, Physics and Mechanics, Jadranska 19, 1000 Ljubljana, Slovenia}
\email{saso.strle@fmf.uni-lj.si}

\thanks{The first author was partially supported by NSF Grant
  DMS-1811111 and NSF grant DMS-1952755.
  The second author was partially supported by ARRS Research program
  P1-0291 and ARRS Grant J1-9104.
	The third author was partially supported by Slovenian Research Agency (ARRS) Research program P1-0288.\\
Math.~Subj.~Class.~2020:  57K40 (primary), 57R40, 57R95, 14J70 (secondary).}

\begin{document}
\begin{abstract} 
What is the simplest smooth simply connected $4$-manifold embedded in $\CP^3$ homologous to a degree $d$ hypersurface $V_d$?  A version of this question associated with Thom asks if $V_d$ has the smallest $b_2$ among all such manifolds.  While this is true for degree at most $4$, we show that for all $d \geq 5$, there is a manifold $M_d$ in this homology class with $b_2(M_d) < b_2(V_d)$.  This contrasts with the Kronheimer-Mrowka solution of the Thom conjecture about surfaces in $\CP^2$, and is similar to results of Freedman for $2n$-manifolds in $\CP^{n+1}$ with $n$ odd and greater than $1$.
\end{abstract}
\maketitle

\section{Introduction}
A conjecture attributed to Ren\'e Thom states that a nonsingular algebraic hypersurface $V_d$ of degree $d$ in $\CP^{n+1}$ is the ``simplest'' representative of its homology class. The notion of complexity of a (real) submanifold $M$ of dimension $2n$ in $\CP^{n+1}$ is motivated by the Lefschetz hyperplane theorem which implies that the homology and homotopy groups of $V_d$ are determined by the ambient manifold $\CP^{n+1}$ below the middle dimension $n$. We are looking for manifolds that closely resemble the behaviour of algebraic hypersurfaces, so the appropriate class of submanifolds within which to look for least complexity representatives of a given codimension-2 homology class in $\CP^{n+1}$ is the class of manifolds $M$ for which the relative (homology and) homotopy groups of the pair $(\CP^{n+1},M)$ are trivial up to the middle dimension. The free parameter is therefore the middle dimensional Betti number, $b_n(M)$.

The conjecture is true for $n=1$, so in $\CP^2$, which was proved by Kronheimer and Mrowka using Seiberg-Witten theory. 

\begin{theorem}[Kronheimer-Mrowka \cite{km-thom}]
A nonsingular algebraic curve of degree $d \in \N$ is a minimal genus smooth surface representing $d[\CP^1]\in H_2(\CP^2;\Z)$;
the first Betti number of such a surface is $d^2-3d+2$.
\end{theorem}

For larger odd $n$ the conjecture is false as was shown by Freedman \cite{freedman}. The idea of the proof is to perform ambient surgery on $V_d$ to reduce the middle Betti number. For technical reasons it is necessary to replace the condition that the relative homotopy groups vanish up to the middle dimension with a stronger condition $\pi_k(C, \partial C)=0$ for $k \le n$, where $C$ is the closed complement of a regular neighborhood of a submanifold $M^{2n} \subset N^{2n+2}$. When this condition holds for the pair $(N,M)$, we say that $M$ is a \emph{taut} submanifold of $N$. Such embeddings were studied by Thomas and Wood \cite{thomas-wood} who in particular showed that $V_d$ are taut and also established lower bounds for $b_n(M)$ for a taut representative $M$ of the class of $V_d$. For $V_d \subset \CP^{2m}$ these lower bounds come from the $G$-signature theorem, are smaller than $b_n(V_d)$ and are almost attained by Freedman's taut manifolds (moreover, he shows they can be realized by rationally taut submanifolds which satisfy the same bounds).

\begin{theorem}[Freedman \cite{freedman}]\label{T:freedman}
For any $m \ge 2$, $V_d$ is not minimal taut in $\CP^{2m}$ for prime $d$, where $d\ne 2,3$ for $m=2$ and $d \ne 2$ for $m=3$.
\end{theorem}

This leaves open the case of $V_d$ in $\CP^{2m+1}$ with $m\geq 1$.  It seems likely that an adaptation of Freedman’s method would establish a result analogous to Theorem \ref{T:freedman} when $m > 1$, so we concentrate on $m=1$.
Thus for any positive integer $d$ we study smooth simply connected 4-dimensional submanifolds $M$ of $\CP^3$ representing the homology class $d[\CP^2]$ which also carry the class of $[\CP^1]$ in $H_2(\CP^3)$. We show that, analogous to the higher dimensional cases, one can find such manifolds $M$ with $b_2(M) < b_2(V_d)$, hence the conjecture does not hold in $\CP^3$. Although we also use ambient surgery we do not need tautness of the embedding due to the special geometric situation in which we perform the construction.

Recall that for a nonsingular $V_d$ the following hold ($\sim$ indicates asymptotic behaviour for large $d$):
\begin{itemize}
\item $V_d$ is simply connected,
\item $b_2(V_d)=d^3-4d^2+6d-2 \sim d^3$, 
\item $\sigma(V_d)=-d(d^2-4)/3 \sim -d^3/3$,
\item $V_d$ is even (spin) for $d$ even, odd for $d$ odd.
\end{itemize}
For small values of $d$ this yields:\\
\begin{center}\begin{tabular}{cccc}
$d$ & $b_2(V_d)$ & $\sigma(V_d)$ & $V_d$\\ \hline
1 & 1 & 1 & \rule{0mm}{5mm}$\CP^2$ \\
2 & 2 & 0 & $S^2 \times S^2$\\
3 & 7 & $-5$ & $\CP^2 \# 6 \overline{\CP^2}$\\
4 & 22 & $-16$ & K3 \\
5 & 53 & $-35$ & quintic
\end{tabular}\end{center}

We show in Proposition \ref{P:basic} that the signature and the parity of its intersection form for a 4-dimensional submanifold of interest in $\CP^3$ are determined by the class it represents and that its $b_2^+$ is at least 1. Inspecting the list above it is then clear that one cannot reduce $b_2$ in any class with $d \le 3$. For $d=4$ the same conclusion follows from the 10/8 Theorem of Furuta \cite{furuta} and in fact from Donaldson's Theorems B and C \cite{donaldson}.

\begin{thm}\label{T:main}
$V_d$ is not minimal in its homology class for $d\ge 5$. There exist simply connected submanifolds $M_d$ of $\CP^3$ homologous to $V_d$ with $b_2(M_d) < b_2(V_d)$. Moreover, for large $d$ we can choose $M_d$ so that $b_2(M_d)$ grows as $3d^3/4$. 
\end{thm}

Using our construction we can reduce $b_2(V_5)$ by 8, so we obtain a $M_5$ with $b_2(M_5)=45$. The smallest $b_2(M_d)$ our method could possibly produce is $\sim d^3/2$ which yields $b_2/|\sigma| \sim 3/2$.  In contrast to Freedman's work, which is restricted to prime degrees, our results work for arbitrary $d\ge 5$. 

As in the work of Freedman and Matsumoto, the proof of the theorem relies on ambient surgery to reduce the second Betti number of the manifold. However, we do not know how to directly implement the approach in \cite{freedman,matsumoto}, so we take a somewhat different route.  
Using results of Baader, Feller, Lewark and Liechti \cite{bfll} we identify a large subgroup of $H_2(V_d)$ on which the intersection pairing is hyperbolic. If these classes were represented by embedded spheres, then they would be candidates for performing surgery on $V_d$ to reduce $b_2$. However, it follows from Donaldson's work (see~\cite[Corollary 6.4.2]{friedman-morgan:book}) that no (non-trivial) homology class of self-intersection $0$ in $V_d$ is represented by an embedded sphere.  On the other hand, Wall~\cite{wall1, wall2} showed that these classes can be represented by spheres after stabilization.  Using this, we can perform ambient surgery to remove a part of the second homology while preserving the characteristic properties of the submanifold.

%
%
\section{Basic properties}

\begin{proposition}\label{P:basic}
Let $d$ be a positive integer and $M \subset \CP^3$ a smooth simply connected 4-dimensional submanifold representing the homology class $d[\CP^2]$ and such that $H_2(M) \to H_2(\CP^3)$ is onto. Then $\sigma(M)=-d(d^2-4)/3$, $M$ is spin iff $d$ is even, and $b_2^+(M) \ge 1$.
\end{proposition}

\proof
By the signature theorem, the signature of $M$ is determined by its first Pontrjagin class. This in turn is determined by the ambient manifold $\CP^3$ and the homology class of $M$, so it agrees with the signature of $V_d$. 

The second Stiefel-Whitney class of the normal bundle of $M$ in $\CP^3$ factors through $H^2(\CP^3;\Z/2)$ and hence is equal to $dx$, where $x$ is the image of the generator of $H^2(\CP^3;\Z/2)$ in $H^2(M;\Z/2)$. Since $\CP^3$ is spin, it follows that $w_2(M)=dx$.

Because the class of $M$ is a positive multiple of $[\CP^2]$ and the inclusion of $M$ into $\CP^3$ induces a surjection on $H_2$, $H_2(M)$ contains a class of positive square.
\endproof

%
%
\section{A model manifold}\label{S:model}

We first choose a smooth algebraic hypersurface $V_d$ of degree $d$ in $\CP^3$ that intersects a 6-ball in a submanifold carrying a large part of the second homology. Recall we are only interested in $d>4$.

\begin{proposition}\label{P:Vd}
$V_d$ can be chosen so that its intersection $F_d$ with a 6-ball $B$ can be isotoped (rel boundary) to the boundary of $B$. Moreover, $F_d$ is the $d$-fold branched cover of the 4-ball branched along a pushed-in Seifert surface $\Sigma_d$ for the $(d-1,d)$ torus knot and $b_2(V_d)=b_2(F_d)+d$.
\end{proposition}

\proof
Let $W_d$ be the singular variety representing the codimension-2 class of multiplicity $d$ in $\CP^3$ given by the equation
$$z_0z_1^{d-1} + z_2^d = z_3^d.$$
Hence $[1:0:0:0] \in W_d$ is the unique singular point (for $d>2$); let $B$ be a small ball about the singularity so that $W_d \cap B$ is the cone on $W_d \cap \partial B$. Clearly $W_d$ is the $d$-fold branched cover of $\CP^2$ with branch set a singular sphere with a unique singular point whose link of singularity is the $(d-1,d)$ torus knot $T_{d-1,d}$. To obtain a smooth representative $V_d$ of the same homology class we choose a nearby nonsingular surface, e.~g.\ the one given by 
$$z_0z_1^{d-1} + z_2^d = \veps z_0^d + z_3^d$$
for a small enough $\veps \ne 0$. In $V_d$ the neighborhood of the singularity $W_d \cap B$ is replaced by the Milnor fibre $F_d$ which can be thought of as the branched cover of $B^4$ with branch set a pushed-in Seifert surface $\Sigma_d$ for $T_{d-1,d}$. An Euler chracteristic computation shows that $b_2(V_d)=b_2(F_d)+d$. Moreover, the Milnor fiber $F_d$ can be isotoped into the boundary sphere of $B$ while fixing its bounday.
\endproof

Next we show there exists a large subgroup of $H_2(F_d)$ (all homology groups from now on have integer coefficients) on which the intersection pairing is hyperbolic. The intersection form of $F_d$ is determined by the Seifert form $\theta_d$ of the Seifert surface $\Sigma_d$. Moreover, $\theta_d$ also determines the linking form $\Theta_d$ on $H_2(F_d) \cong H_1(\Sigma_g) \otimes \Z^{d-1}$ for the embedding of $F_d$ into $\partial B=S^5$; indeed, $\Theta_d=\theta_d \otimes \Lambda_{d-1}$ \cite{durfee-kauffman}, where $\Lambda_k$ is the $k \times k$ matrix of the form
$$\Lambda_k = \left[\begin{matrix} 1 & -1 & 0 & \cdots & 0\\ 0 & 1 & -1 & \cdots &0 \\ \vdots & \vdots & \vdots & \vdots & \vdots \\
0 & \cdots & 0 & 1 & -1 \\ 0 & \cdots & 0 & 0 & 1 \end{matrix}\right].
$$

Even though the smooth slice genus of a torus knot is equal to its genus, the same is not true of its topological (locally flat) slice genus as was first demonstrated by Rudolph \cite{rudolph}. The main tool in the construction is Freedman's result that an Alexander polynomial 1 knot is topologically slice. A systematic study of the topological slice genus of torus knots was conducted by Baader, Feller, Lewark and Liechti \cite{bfll}. They construct subsurfaces of Seifert surfaces whose boundaries are Alexander polynomial 1 knots. We only need the following property of the Seifert form.

\begin{theorem}[\cite{bfll}]
$H_1(\Sigma_d)$ contains a subgroup $G_d$ of rank $2r_d \sim d^2/4$ such that the restriction of the Seifert form $\theta_d$ to $G_d$ is of the following form, consisting of four $r_d \times r_d$ blocks:
$$\left[\begin{matrix} 0 & I + U_d \\  L_d & *\end{matrix}\right],$$
where $U_d$ and $L_d$ are respectively strictly upper- and lower-triangular matrices.
\end{theorem}

\begin{corollary}\label{C:hyperbolic}
The restriction of the Seifert form $\Theta_d$ for $F_d$ to the subgroup $\Ghat_d=G_d \otimes \Z^{d-1}$ of $H_2(F_d)$ has the same form as $\theta_d$ in the previous theorem with the blocks of size $\rhat_d=r_d(d-1)$. Hence the restriction of the intersection form of $F_d$ to $\Ghat_d$ is equivalent to $\oplus{\rhat_d} H$, where $H=\left[\begin{matrix} 0 & 1\\ 1 & 0 \end{matrix}\right]$ denotes the hyperbolic form. The rank of $\Ghat_d$ for large $d$ behaves as $d^3/4$.
\end{corollary}

\proof
Let $(x_i)_i$ be the generators of $\Ghat_d$ corresponding to the first half of the generators for $G_d$ relative to which the Seifert form is given by the matrix in the above theorem and let $(x_i')_i$ be the generators corresponding to the second half. It follows from the structure of the matrix $\Lambda_{d-1}$ that $\Theta_d$ has the same form as $\theta_d$ so it in particular vanishes on the subgroup generated by the $(x_i)_i$. Since the intersection form of $F_d$ is given by $\Theta_d + \Theta_d^\top$, it follows that $x_i \cdot x_j= 0$ and $x_i \cdot x_i' = 1$ for all $i$ and $j$. To make all the other pairings vanish we inductively change the basis elements $(x_i')$ by adding to them linear combinations of $x_j$ for $j \le i$ and $x_j'$ for $j < i$. 
\endproof

%
%
\section{Spherical classes}

In order to reduce the rank of $H_2(F_d)$, we would like to show that some set of generators for the subgroup $\Ghat_d$ of $H_2(F_d)$ can be represented by embedded spheres in $V_d$ and that in fact a regular neighborhood of representatives for a pair of generators giving an $H$ summand as above is diffeomorphic to a punctured $S^2 \times S^2$ and the spheres corresponding to different $H$ summands are disjoint. In general the classes in $\Ghat_d$ may not be represented by embedded spheres (though they are of course spherical) but by Wall's stable diffeomorphism results they are after stabilizing. 

\begin{theorem}[Wall \cite{wall1,wall2}]\label{T:Wall}
Let $M$ and $N$ be simply connected closed 4-manifolds with isomorphic intersection forms. Then the following hold:\\
(1)  for all large enough $\ell > 0$, the stabilized manifolds $M\# \ell (S^2 \times S^2)$ and $N\# \ell (S^2 \times S^2)$ are diffeomorphic;\\
(2) if the intersection form of $M$ is indefinite, any automorphism of the intersection form of $M\# S^2 \times S^2$ is induced by a diffeomorphism.
\end{theorem}

Choose a standard model manifold realizing the intersection form of $V_d$: 
\begin{align*}
M_d & =\frac{b_2+\sigma}2\, S^2 \times S^2 \# |\sigma|\, \overline{\CP^2}\ \ \text{for $d>1$ odd,}\\ 
M_d & =\frac{8b_2+11\sigma}{16}\, S^2 \times S^2 \# \frac{|\sigma|}{16}\, K3\ \  \text{for $d$ even,}
\end{align*} 
where $b_2=b_2(V_d)$ and $\sigma=\sigma(V_d)$. Fix $\ell$ so that $V_d$ and $M_d$ become diffeomorphic after $\ell$ stabilizations. We can realize this stabilization of $V_d$ in $\CP^3$ by internal connected sum of $F_d$ with $\ell$ trivial copies of $S^2 \times S^2 \subset S^5=\partial B$ (each contained in its own 5-disk); denote the stabilized $F_d$ and $V_d$ by $F_d^s$ and $V_d^s$ respectively. We add to $\Ghat_d$ the stabilization classes thus obtaining $\Ghat_d^s \le H_2(V_d^s)$ with $H_2(V_d^s)/\Ghat_d^s\cong H_2(V_d)/\Ghat_d$. 

Denote by $h_d$ the number of $S^2 \times S^2$ summands in $M_d$. Note that $h_d$ for $d$ odd grows as $d^3/3$ whereas for $d$ even as $13d^3/48$, so in any case faster than $\rhat_d \sim d^3/8$. The comparison for small values of $d$ is given in the table below where the data for $r_d$ comes from \cite[Table 1]{bfll}.

\begin{center}\begin{tabular}{ccccccc}
$d$ &\vline & 5 & 6 & 7 & 8 & 9\\ 
$r_d$ &\vline & 1 & 2 & 4 & [5,6] & [6,9]\\
$\rhat_d$ &\vline & 4 & 10 & 24 & [35,42] & [48,72]\\
$h_d$ &\vline & 9 & 9 & 41 & 41 & 113
\end{tabular}\end{center}

We will assume in what follows that $\rhat_d \le h_d$ which is clearly true for large $d$. For those small values of $d$ for which this is not the case we replace $\Ghat_d$ by one of its subgroups satisfying the condition.

\begin{proposition}\label{P:spherical-representatives}
The restriction of the intersection pairing of $V_d^s$ to $\Ghat_d^s$ is equivalent to the sum of hyperbolic forms $H$.  The classes in $\Ghat_d^s$ can be represented by embedded spheres in $V_d^s$ so that for each summand $H$ the corresponding representatives  intersect geometrically once and the spheres corresponding to different $H$ summands are disjoint.
\end{proposition}

\proof
The first claim follows from the construction. By the choice of $\ell$, $V_d^s=V_d\# \ell (S^2 \times S^2)$ is diffeomorphic to the stabilization $M_d^s=M_d\# \ell (S^2 \times S^2)$. We can choose an isomorphism between the intersection pairings of $V_d^s$ and $M_d^s$ that maps the generators of the subgroup $\Ghat^s$ into the generators of the subgroup supported by the sum of $S^2 \times S^2$'s. Since this isomorphism is by Theorem \ref{T:Wall} induced by a diffeomorphism, the second claim follows.
\endproof

%
%
\section{Ambient surgery}\label{S:surgery}

In order to reduce the second Betti number of $V_d^s$ we wish to perform ambient surgery along the spheres guaranteed by Proposition \ref{P:spherical-representatives}. Let $\Sigma_i$ be the spheres representing the first half of the generators of $\Ghat^s$ on which the linking pairing $\Theta_d^s$ vanishes identically (so this collection of spheres contains a representative for one of the generators for each $H$ summand of the restriction of the intersection form to $\Ghat^s$). If $\Sigma_i$ is contained in $F_d^s \subset \partial B$, then it bounds an embedded disk $D_i$ in the 6-ball $B$ and the normal disk bundle of $D_i$ contains an embedded 5-dimensional 3-handle with core $D_i$. The vanishing of the linking pairing guarantees that these handles may be chosen to be disjoint. Since we do not have the control over the action of the diffeomorphism in Wall's stabilization theorem, the spheres might not be contained in $F_d$.  Our main lemma shows that we can arrive at the same conclusion.

\begin{lemma}\label{L:3-disks}
Let $\Sigma_i \subset V_d^s$ be the 2-spheres described above. Then there exist pairwise disjoint embedded 3-disks $D_i \subset \CP^3$ with $D_i \cap V_d^s= \Sigma_i$. Moreover, the disks $D_i$ are not tangent to $V_d^s$.
\end{lemma}

\proof
Denote by $x_i \in H_2(V_d^s)$ the homology class of $\Sigma_i$. Since this class comes from $F_d^s$ it may be represented by an immersed sphere $\Sigma_i^1 \subset F_d^s$ with transverse double points. Then $\Sigma_i^1$ and $\Sigma_i=\Sigma_i^0$ are homotopic in $V_d^s$ (since it is simply connected). According to \cite[Theorem 8.3]{hirsch}, this homotopy may be replaced by a smooth regular homotopy $\varphi_i \colon S^2 \times I \to V_d^s$ (i.e. a homotopy of immersions) if the normal bundle of the immersed sphere $\Sigma_i^1$ is trivial. Since the class $x_i$ has square 0, the latter condition is equivalent to $\Sigma_i^1$ having the same number of positive and negative (transverse) double points. This condition can be satisfied since double points of either sign may be added locally to $\Sigma_i^1$ by replacing a disk with the trace of a homotopy of arcs in $\R^3$ obtained by the sequence of a first Reidemeister move, followed by a crossing change and another first Reidemeister move. (See~\cite[Figure 2]{schwartz:h-cob} for a picture of this process.) We may further assume that the regular homotopy is in general position, so it is a sequence of isotopies, finger moves and Whitney moves \cite[\S 1.6]{freedman-quinn}. The spheres $\Sigma_i^t=\varphi_i(S^2 \times \{t\})$ for $t \in I$ then have transverse double points with the exception of finitely many points; each of these is either the first point of self-intersection for a finger move or the last point of self-intersection for a Whitney move, where the sphere is tangent to itself. Let $\Gamma_i$ be the set of double points of immersed spheres for the map 
$$\Phi_i \colon S^2 \times I \to V_d^s \times I,\quad \Phi_i(x,t)=(\varphi_i(x,t),t).$$
Note that $\Gamma_i$ is the union of properly embedded arcs (with endpoints in $\Sigma_i^1$) and circles; the preimage of $\Gamma_i$ in $S^2 \times I$ consists of two copies of $\Gamma_i$, written $\Gamma_i' \sqcup \Gamma_i''$. 

\begin{figure}[!ht]%
\labellist
\small\hair 2pt
\pinlabel {$\Gamma_i'$} [ ] at 60 185
\pinlabel {$\Gamma_i''$} [ ] at 60 165
\pinlabel {$S^2 \times I$} [ ] at 35 126
\pinlabel {$\gamma_i'$} [ ] at 203 143
\pinlabel {$\gamma_i''$} [ ] at 203 129
\pinlabel {$\Sigma_i^1$} [ ] at 288 92
\pinlabel {$\Sigma_i^0$} [ ] at 14 92
\pinlabel {\rotatebox[origin=c]{270}{$\scaleobj{1.5}{\longrightarrow}$}$\quad\Phi_i$} [ ] at 150 100
\pinlabel {$\Gamma_i$} [ ] at 165 38
\pinlabel {$\gamma_i$} [ ] at 224 31
\endlabellist
\centering
\includegraphics[width=\columnwidth]{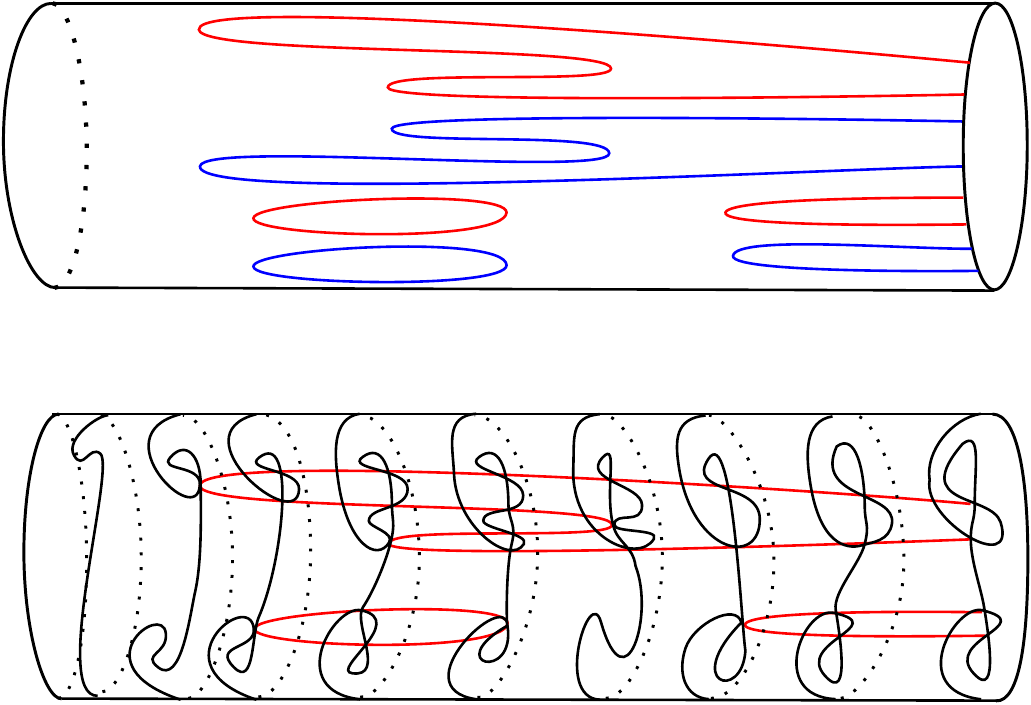}%
\caption{Regular homotopy with separated stages and arcs of double points.}%
\label{F:regular-homotopy}%
\end{figure}
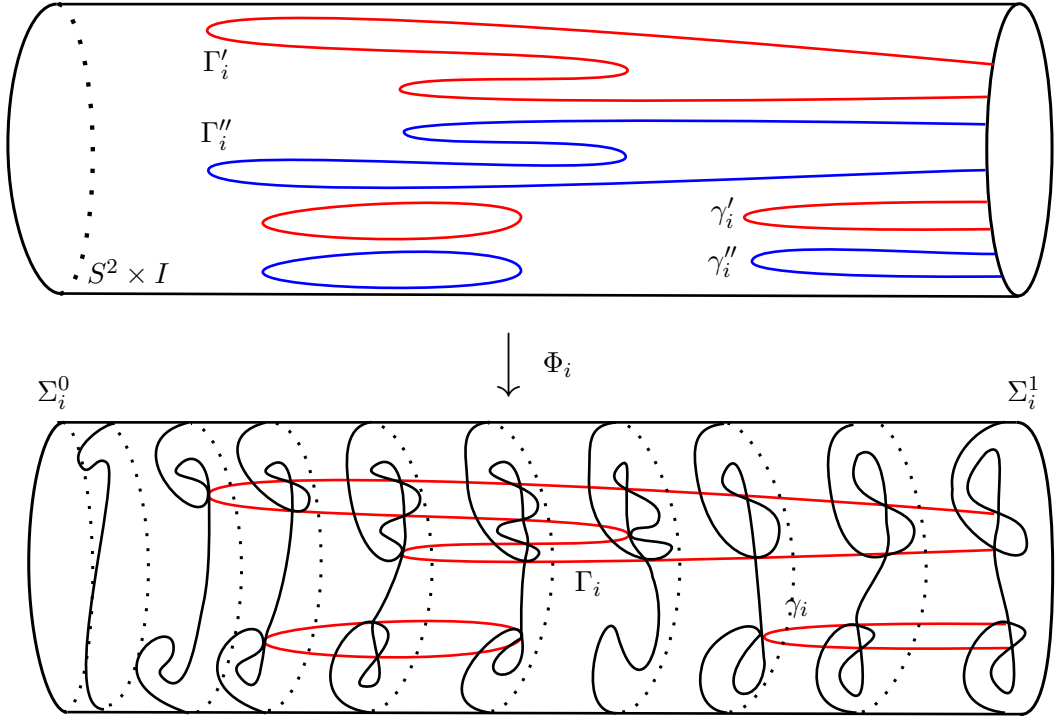

Further, we may assume that the regular homotopies corresponding to different spheres are in general position. This implies that for any time $t$ at most two of the spheres $\Sigma_i^t$ intersect in the same point and this point is not a double point of one of the spheres. The intersections of different spheres are transverse except at tangencies corresponding to finger and Whitney moves. Then $\Delta_{ij}$, the intersection of the images of $\Phi_i$ and $\Phi_j$, is the union of properly embedded arcs (with endpoints in $\Sigma_i^1 \cap \Sigma_j^1$) and circles. The preimages $\Delta_{ij}^i=\Phi_i^{-1}(\Delta_{ij})$ and $ \Delta_{ij}^j=\Phi_j^{-1}(\Delta_{ij}) \subset S^2 \times I$ are two copies of $\Delta_{ij}$.
 
We will push the traces of the homotopies $\varphi_i$ into $\CP^3$ (more precisely, into a tubular neighborhood of $V_d^s$) separating their stages and resolving the (self-)intersections of the immersed spheres thus obtaining disjoint embedded annuli $W_i$ connecting $\Sigma_i^0$ to a push-off of $\Sigma_i^1$. Equip the normal bundle $\nu$ of $V_d^s \subset \CP^3$ with a Riemannian metric; we may assume that the metric over $F_d^s$ is induced from the metric on the 6-ball $B$, which we identify with the ball of radius 2 in $\R^6$. By rescaling the metric (by a constant factor) we may assume that the unit disk bundle of $\nu$ is identified with a tubular neighborhood of $V_d^s$; we will use this identification implicitly in what follows.

We fix a trivialization of the pull-back of $\nu$ via $\varphi_i$. This bundle is the pull-back of the trivial bundle over $S^2 \times \{1\}$, and there is a particular choice of trivialization over this sphere given by the vector fields $E_1^i$, the pull-back of the inner normal to $S^5$ in $B^6$, and $E_2^i$, the pull-back of the normal vector field to $F_d^s$ in $S^5$. Choosing a trivialization of the bundle over $S^2 \times I$, we extend $(E_1^i,E_2^i)$ to orthonormal trivializing sections $(E_1^i,E_2^i)$ of the whole bundle. Let $\lambda \colon [0,1] \to [0,1]$ be a smooth increasing surjective function that is constant in some neighborhoods of the endpoints. Then $\psi_i \colon S^2 \times I \to \CP^3$, given by
$$(x,t) \mapsto t E_1^i (\varphi_i(x,\lambda(t))),$$
is an embedding of the image of $\Phi_i$ in $\CP^3$ (with collars added at each end). Denote the image of $\psi_i$ by $Z_i$. Note that $\varphi_i$ factors through $Z_i$, where $Z_i$ maps to $V_d^s$ by the projection. In particular, the pull-back of $\nu$ via $\varphi_i$ factors through its pull-back to $Z_i$. So for any component $\gamma_i$ of $\Gamma_i \subset Z_i$ we may identify the pull-back of $\nu$ to $\gamma_i$ with its pull-back to either component $\gamma_i'\subset \Gamma_i'$ or $\gamma_i''\subset \Gamma_i''$ of its preimage.

In order to get embedded annuli $W_i$ we first need to remove the double points of immersed spheres. Note that over any corresponding pair of components $\gamma_i' \sqcup \gamma_i''$ in $\Gamma_i' \sqcup \Gamma_i''\subset S^2 \times I$ that map to $\gamma_i \subset \Gamma_i \subset Z_i$, the two trivializations of the pull-back of the normal bundle $\nu$ to $\gamma_i$ determined by $(E_1^i,E_2^i)$ restricted to either $\gamma_i'$ or $\gamma_i''$ are homotopic as any such component is null-homotopic in $S^2 \times I$. We now change $E_1^i$ over $\gamma_i''$ to agree with the restriction  of $E_1^i$ to $\gamma_i'$ rotated by a small angle $\delta>0$ in the direction of $E_2^i|\gamma_i'$. Choose small pairwise disjoint compact regular neighborhoods $K_i=K_i' \sqcup K_i''$ for $\Gamma_i' \sqcup \Gamma_i''$ and $L_{ij}^i \sqcup L_{ij}^j$ for $\Delta_{ij}^i \sqcup \Delta_{ij}^j$. Using $(E_1^i,E_2^i)|K_i'$ to trivialize the normal bundle over $K_i$, we choose the fibrewise universal cover of the corresponding circle bundle in which $E_1^i|K_i'$ corresponds to the zero and $E_2^i|K_i'$ to a positive angle. Then the lift of $E_1^i|K_i''$ may be smoothly spliced with the constant section $\delta$ and then pushed down into the circle bundle to give the new section $E_1^i|K_i''$; then rotating $E_2^i|K_i''$ appropriately we obtain an othonormal frame. In fact, when $\gamma_i$ is an arc, we complete $\gamma_i''$ to a circle by adding to it an arc in $S^2 \times \{1\}$. Since over $S^2 \times \{1\}$ the section $E_1^i$ is determined as the inner normal to the boundary of $B$, we choose the lift of $E_1^i|K_i''$ over this arc to be zero. This shows that the double points of $\Sigma_i^1$ are removed by a small homotopy inside the ball $B$.

\begin{figure}[!ht]%
\includegraphics[width=0.5\textwidth]{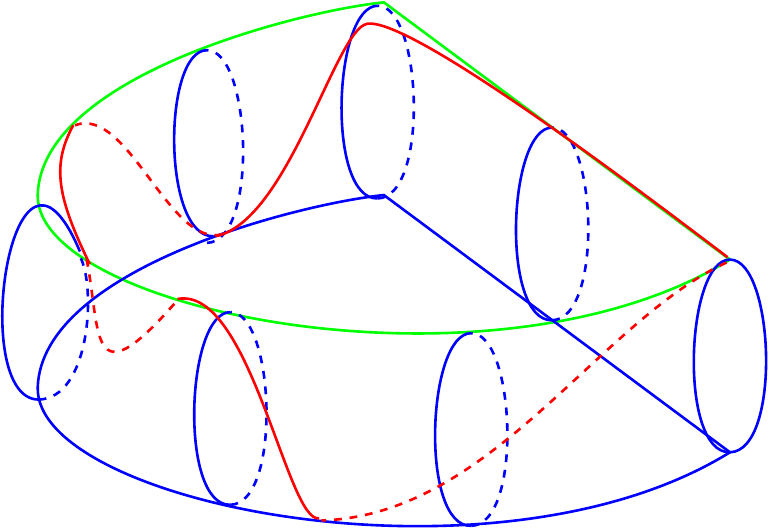}%
\caption{Comparison of the two trivializations of the pull-back of $\nu$ to $\gamma_i$, given by the sections $E_1^i|\gamma_i'$ (green) and $E_1^i|\gamma_i''$ (red). In the picture, $\gamma_i$ is assumed to be an arc and is completed to a circle by an arc over which the two sections agree.}%
\label{F:trivializations}%
\end{figure}

To remove the intersections between different spheres we just repeat the same procedure with any $\Delta_{ij}$, where we assume $i<j$, by changing the section $E_1^j$ over $L_{ij}^j$. Denote the resulting embeddings obtained in this way from $\psi_i$ by $\Psi_i \colon S^2 \times I \to \CP^3$; then one boundary component of $W_i = \Psi_i(S^2 \times I)$ is $\Sigma_i$ and the other is the push-in of $\Sigma_i^1$ which we denote $\wSig_i=\Psi_i(S^2 \times \{1\})$.
Note that the spheres $\wSig_i$ are essentially contained in a 5-sphere $S$ concentric with the boundary of the ball $B$; they only deviate from $S$ in small neighborhoods of double points and intersection points between different spheres (more precisely, over the images of $K_i'' \cap \wSig_i$, and over $L_{ij}^j \cap \wSig_j$ for $i < j$). But as noted above, the removal of intersection points in $\Sigma_i^1$ is realized by a small homotopy. Hence the projection of $\wSig_i$ into $S$ along the normal vector field is a diffeomorphism and we may and will assume that $\wSig_i$ is contained in $S$. Then $\wSig_i$ bounds a properly immersed 3-disk $D_i$ in the ball $B'$ bounded by $S$. Assuming $D_i$ is in general position, it may have transverse self-intersections, but pairs of double points in $D_i$ of opposite sign can be canceled using the Whitney trick. Note that the number of double points of either sign may be increased by adding kinks (analogous to first Reidemeister move) into $\wSig_i$. Thus we may assume that $D_i$ is embedded.

Since $\Theta_d^s(\Sigma_i^1,\Sigma_j^1)=0$ and $\wSig_i$ is just a push-in of $\Sigma_i^1$ into a concentric sphere, the linking number of $\wSig_i$ and $\wSig_j$ (where the latter can be considered as a perturbation of the normal push-off of $\Sigma_j^1$) is trivial. Hence the intersection number of $D_i$ and $D_j$ is zero for all $i \ne j$, so we may assume that they are geometrically disjoint (by using the Whitney trick).
\endproof

\begin{theorem}
The homology class of $V_d$ in $\CP^3$ is represented by a simply-connected manifold $N_d$ with $H_2(N_d) \cong H_2(V_d)/\Ghat$.
\end{theorem}

\proof
We first show that the disk $D_i$ may be thickened to a 5-dimensional 3-handle $h_i$ in $\CP^3$ with the attaching region contained in $V_d^s$ and whose attaching sphere is equal to $\Sigma_i$.
The normal bundle of $\Sigma_i$ in $V_d^s$ is trivial, so its normal disk bundle in $\CP^3$ admits a splitting $\Sigma_i \times B^2 \times B^2$, where the first $B^2$ corresponds to the normal directions in $V_d^s$, and the second corresponds to the restriction of the normal bundle $\nu$ of $V_d^s \subset \CP^3$ to $\Sigma_i$. The latter is trivialized by $(E_i^1,E_i^2)$ and the normal disk bundle of $D_i$ over $\Sigma_i$ is given by $\Sigma_i \times B^2 \times B^1 E_i^2$.  This trivialization extends over the normal bundle of $D_i$ in $\CP^3$ since $\pi_2(\GL_3\R)$ is trivial. The required handle $h_i$ is $D_i \times B^2 \times 0$.

Using the handles $h_i$ we perform ambient surgery on $V_d^s$ along the $\Sigma_i$ to obtain a manifold $N_d$, homologous to $V_d^s$ and hence to $V_d$. Clearly $H_2(N_d)$ is isomorphic to $H_2(V_d^s)/\Ghat^s\cong H_2(V_d)/\Ghat$ since surgery on $\Sigma_i$ kills also its dual class. 

That $N_d$ is simply connected follows since the fundamental group of the complement of $\Sigma_i$ is normally generated by its meridian which is trival in $N_d$, because the dual class to $\Sigma_i$ is also represented by a sphere.
\endproof 

The final question to address is whether the manifolds $N_d$ are taut in $\CP^3$. We show below that $\pi_k(\CP^3,N_d)$ is trivial for $k \le 2$. In fact, it also follows by general position arguments that $\pi_1(C,\partial C)$ is trivial, where $C$ is the closed complement of a tubular neighborhood of $N_d$. 
\begin{proposition}
The pair $(\CP^3,N_d)$ is $2$-connected.
\end{proposition}

\proof
Since $\CP^3$ and $N_d$ are simply-connected, we only need to verify that the inclusion induced homomorphism is surjective on $\pi_2$ or equivalently on $H_2$. Since $V_d$ is taut, so is the stabilized manifold $V_d^s$ (by the argument as in the previous sentence). Hence the generator $x \in H_2(\CP^3)$ is the image of an element $\widetilde x \in H_2(V_d^s)$. Suppose now we do the surgery on a sphere $\Sigma_i$ representing the class $x_i \in H_2(V_d^s)$. Let $y_i \in H_2(V_d^s)$ be the homology class of the dual sphere to $\Sigma_i$. Then $x_i$ has trivial algebraic intersection with the class
$$\widetilde x' := \widetilde x -(\widetilde x \cdot x_i)y_i.$$
Since $y_i$ is supported by the Milnor fibre $F_d \subset B$, it maps to the trivial class in $H_2(\CP^3)$, so the surgery preserves surjectivity. 
\endproof
In general, the ambient surgery construction may destroy tautness, but based on a theorem of Quinn~\cite{quinn:canonical} as quoted by Freedman \cite[Theorem 2.5]{freedman},  it seems that one can re-embed $N_d$ as a taut submanifold.\\[1ex]
{\bf Remarks on the proof.} The overall strategy used to prove Theorem~\ref{T:main} is similar to that in the work of Freedman and Matsumoto, but with an important difference. In our work and also in~\cite{freedman,matsumoto} an algebraic form on a subspace of the middle homology is an obstruction to doing ambient surgery.  (For us it is essentially the Seifert form, whereas~\cite{freedman,matsumoto} use a Wall-type form denoted $(\lambda, \mu)$.)  However, the technique in~\cite{freedman,matsumoto}, applied in our setting, would be to immerse a $3$-disk with boundary on $V_d$, and use the vanishing of $(\lambda, \mu)$ to push the singularities out to the boundary $2$-spheres; these would in principle be removed by an application of the Whitney trick. Since the Whitney trick does not apply in dimension $4$, we modified the procedure to get embedded $2$-spheres (after stabilization) and then remove the singularities of the $3$-handles by the Whitney trick in dimension $6$.

\end{document}